\documentclass[12pt]{article}
\usepackage{amsmath,amssymb,amsthm,mathtools,url}

\newtheorem{prob}{Problem}

\newcommand{\zp}{\mathbb{Z}/p\mathbb{Z}}
\def\({\Big[}
\def\){\Big]}
\def\mand{\,\,\,\land\,\,\,}
\def\id{1'}
\def\P{\mathrm{Pr}}

\title{Finite representations for two small relation algebras}

\author{Jeremy F. Alm\\Department of Mathematics\\Lamar University \\ Beaumont, TX 77710\\
\texttt{alm.academic@gmail.com}\\
\and Roger D. Maddux\\Department of Mathematics\\Iowa State University\\
\texttt{maddux@iastate.edu} }
\date{September 2017}

\begin{document}

\maketitle

\begin{abstract}
In this note, we give two different proofs that relation algebra $52_{65}$ is representable over a finite set.  The first is probabilistic, and uses Johnson schemes.  The second is an explicit group representation over $ (\mathbb{Z}/2\mathbb{Z})^{10}$. We also give a finite representation of $59_{65}$ over $\mathbb{Z}/113\mathbb{Z}$ using a technique due to Comer.
\end{abstract}

\section{Introduction}

There are 65 finite integral symmetric relation algebras with exactly three diversity atoms. In Maddux's book  \cite{Madd}  they are numbered $1_{65}$ through $65_{65}$.  For each such algebra, it is known whether it is representable or not.  For those that are representable, whether representations are possible over finite sets is trickier, and it is not known for every such algebra.  In this note we show that $52_{65}$ and $59_{65}$ are finitely representable.

\section{Relation algebra $59_{65}$}

Relation algebra $52_{65}$ has atoms $\id$, $a$, $b$, $c$, and diversity cycles $aaa$,
$bbb$, $acc$, $aab$, $aac$, $bcc$, $abc$, but not $ccc$, $abb$, or $cbb$. (The identity
cycles are $1'1'1'$, $1'aa$, $1'bb$, and $1'cc$.)   In any representation, $b\cup\id$ is
therefore an equivalence relation on $U$. We prove the existence of a representation with exactly
three equivalence classes.  The argument is similar to the one given in \cite{AMM}. Choose $n\in\omega$ large enough for what follows. Let $U$ be
the $n$-element subsets of a set with $3n-4$ elements. Randomly partition $U$ into three
sets $S_0$, $S_1$, $S_2$. Let
\begin{align*}
	\id&=\{(x,x):x\in U\}
\\	b&=\bigcup_{i=0}^2\{(x,y):x,y\in S_i,\,x\neq y\}
\\	a&=\{(x,y)\in U^2: |x\cap y|\geq2\}
\\	c&=\{(x,y)\in U^2: |x\cap y|\leq1\}
\end{align*}

Let $\Phi_{xyz}^{ijk}$ stand for the open formula 
\[
    (x,y)\in i\mand ijk\text{ is a cycle }\mand
	\big( (x,z)\notin j\lor (z,y)\notin k \big),
\]
i.e., $z$ fails to witness the need $ijk$ for the edge $xy$ colored $i$. 
The probability (which is $<1$ for all large enough $n$) that this is NOT a representation is
\begin{align*}
  &	\P\Bigg( \(\exists x,y\in U\) \(\exists i,j,k\in\{1',a,b,c\}\) 
	\(\forall z\in U\) \ \ \Phi_{xyz}^{ijk} \Bigg) \\
&\leq	\sum_{x,y\in U} \P\Bigg( \(\exists i,j,k\in\{1',a,b,c\}\) 
	\(\forall z\in U\)\Bigg)  \ \ \Phi_{xyz}^{ijk} \Bigg) \\
&\leq	\sum_{\mathclap{\substack{x,y\in U\\ i,j,k\in\{1',a,b,c\}}}}
	\P\Bigg(\(\forall z\in U\)  \ \  \Phi_{xyz}^{ijk}\Bigg) \\
&=\sum_{\mathclap{\substack{x,y\in U\\ i,j,k\in\{1',a,b,c\}}}} \ \ \prod_{z\in U}
	\P\Bigg( \Phi_{xyz}^{ijk} \Bigg)\\
\end{align*}

The situation with the smallest number of witnesses is $(x,y)\in a$ with $|x\cap y|=2$, cycle
$acc$. A witness is some $z\in U$ with $(x,z),(z,y)\in c$, so $z$ is not in the equivalences
class(es) of $x,y$. The number of such $z$ is smallest if $x,y$ are in distinct classes.  There are
$(n-2)^2$ possible $z$'s, and the probability that all of them are in the two equivalence classes
of $x,y$ is $(2/3)^{(n-2)^2}$---if we use three classes \emph{of the same size}! So let's assume
that $|U|$ is divisible by 3, and pick a random partition of $U$ that has three sets of the same
size---each is a \emph{third of $U$}. The number of such partitions is
$\frac12{|U|\choose|U|/3}{2|U|/3\choose|U|/3}$ where $|U|={3n - 4 \choose n}$. Then we have 
\begin{align*}
  & \sum_{\mathclap{\substack{x,y\in U\\ i,j,k\in\{1',a,b,c\}}}} \ \ \prod_{z\in U}
	\P\Bigg( \Phi_{xyz}^{ijk} \Bigg)\\
&\leq	\sum_{\mathclap{\substack{x,y\in U\\ i,j,k\in\{1',a,b,c\}}}} \ \ (2/3)^{(n-2)^2}	\\
&\leq	|U|^2\cdot4^3\cdot(2/3)^{(n-2)^2}	\\
&=	{3n - 4 \choose n}^2\cdot4^3\cdot(1-1/3)^{(n-2)^2} < 1	
\end{align*}
if $n$ is large enough. A little Python script informs us that $n=13$ suffices, giving   ${3n - 4 \choose n}= 1476337800$.

Now, clearly, 1.4 billion points is more than necessary. Let's try abelian groups.  If we are to partition an abelian group $G=\{0\}\cup A \cup B \cup C$ (with $A$ being the image of $a$, etc) then $B\cup \{0\}$ is a subgroup and $C$ is a maximal sum-free set, i.e., $C+C=G\setminus C$. Consider $G=(\mathbb{Z}/2\mathbb{Z})^{10}$, and consider the elements as bitstrings. Let 
\[
    X=\{x\in (\mathbb{Z}/2\mathbb{Z})^{10} : x \text{ has between one and six 1s}\}
\]
and  
\[
    C=\{x\in (\mathbb{Z}/2\mathbb{Z})^{10} : x \text{ has more than six 1s}\}.
\]
It is not hard to check that $X+X=G$, $X+C=G\setminus\{0\},$ and $C+C=G\setminus C$.  

Now we seek a subgroup $H\subseteq X\cup \{0\}$, and let $B = H\setminus\{0\}$ and $A=X\setminus B$.  If $H$ is large enough, we should have $A+B=B+C=A\cup C$.  (We should also check $A+A=G$ and $A+C=A\cup B \cup C$, but those  will almost certainly hold simply by cardinality considerations.)

The following is a list of the elements of $B$, the nonzero elements of a subgroup $H$ that works.  This relatively large subgroup was found by computer-assisted trial and error. $H$ induces 16 equivalence classes, each one a $B$-colored clique of order 64.

\begin{quote}
{\footnotesize 
0110000000, 
1000000010, 
0000000101, 
1000000111, 
0100001000, 
1100001010,

0100001101, 
1100001111, 
0000010001, 
1000010011, 
0000010100, 
1000010110,

0010011001, 
0100011001, 
1010011011, 
0100011100, 
1100011110, 
0110000101, 

1000100000, 
0000100010, 
1100011011, 
1000100101, 
0000100111, 
1100101000, 

0010011100, 
0100101010, 
1110000111, 
1100101101, 
0100101111, 
1000110001, 

0000110011, 
1000110100, 
1010011110, 
0000110110, 
1100111001, 
0100111011, 

1100111100, 
1010001010, 
0100111110, 
1110100000, 
0010101111, 
0110110110, 

1110110100, 
0110110011, 
1110010110, 
0110100010, 
0010001101, 
1010111001, 

1010001111, 
1110000010, 
1110100101, 
1010101101, 
0010111011, 
1110110001, 

0110010001, 
1010111100, 
0110100111, 
1010101000, 
1110010011, 
0010111110,

0110010100, 
0010001000, 
0010101010 }
\end{quote}


\section{Relation algebra $59_{65}$}

Relation algebra $59_{65}$ has atoms $\id$, $a$, $b$, $c$, and diversity cycles $aaa$,
 $acc$, $aab$, $aac$, $bcc$, $abc$, and  $ccc$, but not $bbb$ or $cbb$.   For $59_{65}$, we use a construction borrowed from Comer \cite{Comer, AlmYlv}. 

Let $p=113$, and fix $m=8$. For $0\leq i <m$, let 
\[
    X_i = \{g^{\alpha m + i} : 0\leq \alpha < (p-1)/m \},
\]
where $g=3$ is the smallest primitive root modulo $p$.  $X_0$ is the (unique) multiplicative subgroup of $(\zp)^\times$ of index 8, and the other $X_i$s  are its cosets. The $X_i$s form atoms of a relation algebra with forbidden cycles $ [X_{i}, X_{i}, X_{i}]$,  $[X_{i}, X_{i}, X_{i+6}]$, and $[X_{i}, X_{i}, X_{i+7}]$.  The cycle structure does not depend on the choice of generator $g$, but the indexing of cosets does. 

Let $A = X_1 \cup X_2 \cup X_3 \cup X_4 \cup X_5 $, and let $B = X_0$ and  $C = X_6 \cup X_7$. 

Then we have the following:

\begin{itemize}
    \item $A+A=\zp$
    \item $A+B=(\zp)^\times$
    \item $A+C=(\zp)^\times$
    \item $B+B=\{0\}\cup A $
    \item $B+C=A\cup C $
    \item $C+C=\zp$
\end{itemize}

\section{Concluding remarks}

Now we may cross $52_{65}$ and $59_{65}$ off of the list of integral relation algebras with three symmetric diversity atoms that are representable, but not known to be representable on a finite set. The four remaining algebras (to the best of our knowledge) are $30_{65}$, $33_{65}$, $34_{65}$, and $56_{65}$.  

\begin{prob}
Are any of relation algebras  $30_{65}$, $33_{65}$, $34_{65}$, and $56_{65}$ representable over a finite set?  Over a finite group?

\end{prob}

\begin{prob}
Is $52_{65}$ representable over $ (\mathbb{Z}/2\mathbb{Z})^{7}$ by a similar method? (It is important here that the exponent is equivalent to $1\pmod{3}$.)
\end{prob}

\begin{prob}
Is $52_{65}$ representable over $ (\mathbb{Z}/2\mathbb{Z})^{3k+1}$ for arbitrarily large $k$?  We were able to find a representation over $ (\mathbb{Z}/2\mathbb{Z})^{13}$, again by trial and error.  The difficulty is in constructing the subgroup $H$. 
\end{prob}



\end{document}